\documentclass[11pt,reqno]{amsart}
\usepackage{graphicx}
\pdfoutput=1
\usepackage{amsfonts,amsmath,amssymb}
\usepackage{hyperref}
\usepackage{paralist}
\usepackage[linesnumbered,ruled,vlined,norelsize]{algorithm2e}
\usepackage{graphicx}

\newcommand{\cx}{{\mathbb{C}}}

\newcommand {\Q}{\mathcal Q}

\usepackage{thmtools}
\declaretheoremstyle[bodyfont=\normalfont]{noncursive}
\declaretheorem{theorem}
\declaretheorem[numberwithin=section]{lemma}

\declaretheorem[numberlike=lemma]{proposition}

\declaretheorem[style=noncursive,numberlike=lemma]{definition}

\declaretheorem[style=noncursive,numberlike=lemma]{remark}

\newcommand{\sss}{\mathcal{S}}
\newcommand{\CC}[1]{\mathbb{C}^{#1}}
\newcommand{\CP}[1]{\mathbb{CP}^{#1}}

\numberwithin{equation}{section}

\title[Sphericity of a real hypersurface]{Sphericity of a real hypersurface via projective geometry}

\author {Ilya Kossovskiy}
\address{\parbox{0.8\linewidth}{%
        Department of Mathematics, Federal University of Santa Catharina/\\ %
        Department of Mathematics and Statistics, Masaryk University, Brno}
    }
\email{ilyakos@gmail.com}

\begin{document}

\maketitle

\date{\today}

\begin{abstract}
In this work, we obtain an unexpected geometric characterization of  sphericity of a real-analytic Levi-nondegenerate hypersurface $M\subset\CC{2}$. We prove that $M$ is spherical if and only if its Segre\,(-Webster) varieties satisfy an elementary combinatorial property, identical to a property of straight lines on the plane and known in Projective Geometry as the {\em Desargues Theorem}.  
\end{abstract}

\tableofcontents

\section{Introduction}

Let $M\subset\CC{2}$ be a real-analytic Levi-nondegenerate hypersurface. The celebrated theory due to Chern and Moser (see also earlier work of Poincar\'e \cite{poincare} and Cartan \cite{cartan}) gives a beautiful exposition of the fact that only a very rare such $M$ is locally biholomorphically equivalent to the 3-dimensional sphere $S^3\subset\CC{2}$. It also demonstrates the exceptional role of the sphere $S^3$ as the {\em model} for {all} Levi-nondegenerate hypersurfaces in $\CC{2}$.   One of the central questions arising in connection with that is the problem of an {\em effective} identification of those Levi-nondegenerate hypersurfaces  which are locally biholomorphically equivalent to $S^3$ (so-called {\em spherical hypersurfaces}). Chern and Moser provide a certain answer in this direction. First, the differential-geometric construction of Chern (see also \cite{cartan}) gives an algorithm for computing the {\em CR-curvature} of a Levi-nondegenerate hypersurface, vanishing of which is equivalent to the  sphericity. (We note that the local sphericity of a real-analytic Levi-nondegenerate  hypersurface is equivalent to its global sphericity, see, e.g., \cite{pinchuk}). On the other hand, Moser provides a (convergent) normal form construction for a real-analytic Levi-nondegenerate hypersurface. In this language, the sphericity reads as  vanishing of all the ``resonant'' terms $\Phi_{kl}(u)$ Moser's distinguished coordinates (see the formula \eqref{nf} in Section 4 below). 

One of the remaining problems in Chern-Moser's theory is to find an {\em elementary geometric} characterization of sphericity of a real hypersurface.  That is, we are searching for elementary geometric properties of a real hypersurface $M\subset\CC{2}$,  guaranteeing its sphericity  on the one hand, and {not requiring} a transfer to any special coordinates on the other hand. Somewhat surprisingly, {\em no} results in this direction exist till present, and it was commonly expected that the sphericity can {\em not} be characterized any simpler than via Chern's tensor invariants or Moser's normal form.

An elementary geometric characterization under discussion is not only interesting by itself, but is also motivated by the fact that both Chern's curvature and Moser's normal form constructions are usually difficult to carry out for a concrete Levi-nondegenerate hypersurface due to their very high computational complexity. Checking the sphericity even of very simple  hypersurfaces (e.g., real-algebraic ones) remains a very difficult task in general   (the complexity of Chern's approach can be seen, for example,  from the respective computation for real ellipsoids  \cite{huangumbilic,websterduke}).

\smallskip

The main goal of this paper is to provide such a, very unexpected, geometric characterization of sphericity. The associated geometric object that we use is a useful tool due to Segre and Webster which we address as {\em the Segre-Webster family} of a real hypersurface (note that the latter object is often  addressed in the literature as merely {\em the Segre family}; we refer to \autoref{web} below for the history and terminology, and to Section 2 for details of the concept and properties). 

\smallskip

Our main result below further demonstrates the exceptional role of Segre-Webster families for CR-geometry. Informally speaking, the result says that  a real-analytic Levi-nondegenerate hypersurface $M\subset\CC{2}$ is spherical near a point $O\in M$ if and  only if its Segre-Webster family near $O$ (which is a $2$-parameter family of planar complex curves) satisfies one of the configuration theorems of Projective Geometry: the {\em Desargues Theorem}. For convenience of the reader, we provide below a picture illustrating the Desargues theorem. 

\bigskip

\includegraphics[scale=0.42]{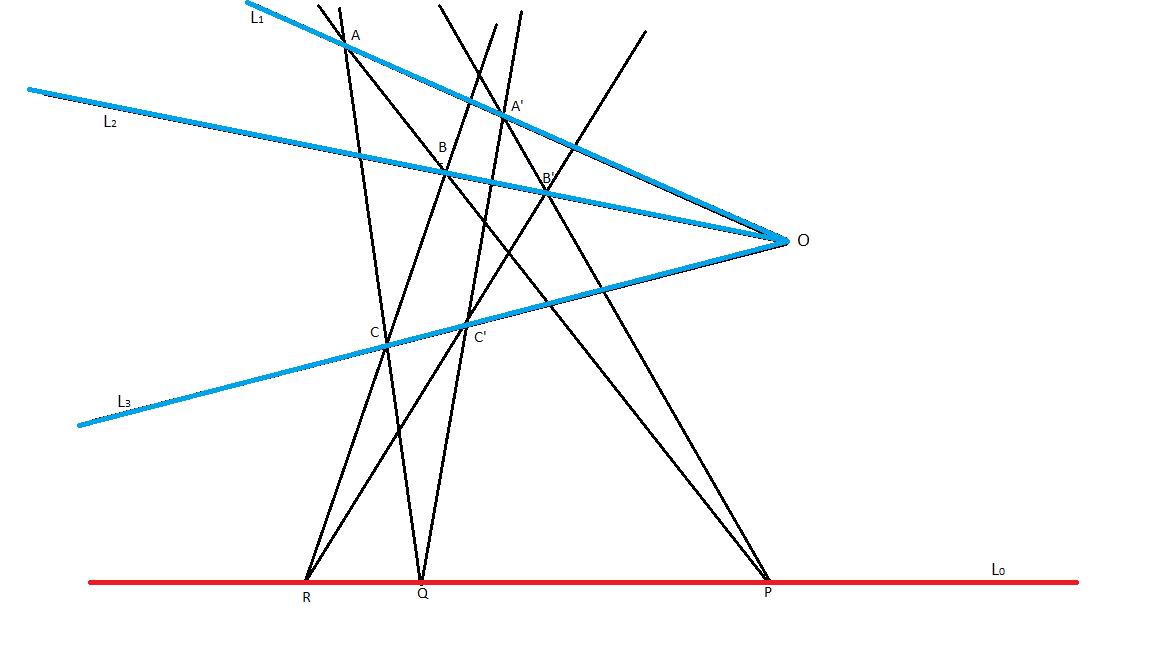}

\bigskip

The lines on the picture are arbitrary real or complex lines in $\mathbb{RP}^2$ or $\CP{2}$ respectively. The assertion of the theorem is that the three points $$P,Q,R,$$ obtained from the three distinct lines $L_1,L_2,L_3$ with $L_1\cap L_2\cap L_3=\{O\}$ and three pairs of  points $A,A'\in L_1,\,B,B'\in L_2,\,C,C'\in L_3$ as shown on the picture,

\begin{center}{\em must lie on the same line $L_0$.}\end{center}

\noindent Here the seven points $A,A',B,B',C,C',O$ are assumed to be pairwise distinct. 

\begin{definition}
In what follows, we call the latter property of a configuration of seven points {\em the Desargues property}.
\end{definition}


Let us now consider a more general (local) setting, where the family of lines on the plane is replaced by a (local!) $2$-parameter holomorphic family of complex curves, behaving similarly to the family of lines on the plane. More precisely, we consider a  family of complex curves $\sss=\{L_{c}\}$ parameterized by a parameter $c$ running an open set $\Omega\subset\CC{2}$,  such that: 

\bigskip

\noindent (i) each of the curves $L_{c}$ is defined in a polydisc $U\subset\CC{2}$ centered at a point $O$;

\bigskip

\noindent (ii) the holomorphic parameter $c$ lies in a polydisc  $V\subset\CC{2}$ centered at a point $c_0$, $L_{c_0}\ni O$ and the parameterization map $c\mapsto L_c$ is holomorphic and one-to-one in $V$;

\bigskip

\noindent (iii) for all $p\in U$ and all slopes $l$ sufficiently close to the slope $l_0$ of $L_{c_0}$ at $O$, there exists one and only one curve $L_c$ from $\sss$, passing through $p$ and having the slope $l$ at $p$.   

\bigskip

In other words, the lifting of the $2$-parameter family $\sss$ to the bundle of $1$-jets of planar complex curves forms a foliation near a point in the $1$-jet bundle corresponding to $O,l_0$.

\begin{definition}
In what follows we call every family satisfying (i)-(iii) {\em a transverse family}. 
\end{definition}

As a well known fact (e.g., \cite{ber},\cite{DiPi}), 

\bigskip

{\em the Segre-Webster family of a real-analytic Levi-nondegenerate hypersurface $M\subset\CC{2}$ near every point $O\in M$ is a transverse family}. 

\bigskip

As follows from the definition of a transverse family, after possibly shrinking the polydiscs $U,V$, we have the following property: through every two points $A,B\in U$ there exists at most one curve of the family passing through $A,B$. Thus, we can correctly define the unique intersection point $L_1\cap L_2$ of two curves $L_1,L_2$ in a transverse family $\sss$, as well as the unique curve $L=AB$ passing through two given points $A,B\in U$, if the latter objects exist. 

\begin{definition} We call a configuration of seven distinct points $O,A,B,C,A',B',C'$ in $U$ {\em an admissible configuration}, if each of the triples of points $(O,A,A'),(O,B,B'),(O,C,C')$ belongs to some curve  in $\sss$, all the intersections $P=AA'\cap BB'$, $Q=AA'\cap CC'$, and $R=BB'\cap CC'$ are non-empty, and all the three curves $PQ,PR$, and $PR$ are non-empty. 
\end{definition}

It then makes a complete sense, replacing  ``lines'' by ``cuves of a transverse family $\sss$'', to ask whether an admissible configuration has the Desargues property in the above sense (that is, whether $P,Q,R$ lie in the same curve in the family $\sss$).  We now do the following simple observation: if a transverse family of curves is locally mappable onto the family of straight lines, then (after shrinking the basic neighborhood) for  all admissible configurations for $\sss$ the Desargues property must hold (as it holds for the transformed family).  In particular, if $\sss$ is the Segre-Webster family of a real hypersurface $M$, then in the new coordinates the Segre-Webster varieties are straight lines and, as a well known fact, the image of $M$ lies in a real hyperquadric in $\CP{2}$. Hence  $M$ is spherical, and we conclude that 

\smallskip

{\em the Desargues property is necessary for the sphericity}.   

\smallskip

Our main result says that, in fact, the Desargues property  is also sufficient here.

\begin{theorem}\label{theor1}
Let  $M\subset\CC{2}$ be a real-analytic Levi-nondegenerate hypersurface, and $p\in M$. Then the germ $(M,p)$ is spherical if and only if the Segre-Webster family of $M$ near $p$ has the Desargues property for all admissible configurations.
\end{theorem}

\begin{remark}\label{web} Segre-Webster families (normally addressed in the literature as merely {\em Segre families}) were discovered in the work \cite{segre} of Segre, and revisited more recently by Webster \cite{webster}. Webster used them in \cite{webster} for proving several celebrated theorems in CR-geometry concerning mappings of algebraic hypersurfaces. Since the work \cite{webster}, Segre families is a key tool for  modern CR-geometry (see, e.g., \cite{ber},\cite{DiPi} and references therein).  During the conference in Madison, Wisconsin  in March 2015 dedicated to the 70th anniversary of  Webster, we agreed with my colleagues that, taking into account the exceptional role of Webster in discovering and pursuing Segre families as one of the central objects in CR-geometry,  it would be fair to address Segre families as {\em Segre-Webster families}. Thus, throughout the paper we use the latter term only.
\end{remark}  


The paper is organized as follows. In Section 2 we give a background information on Segre-Webster families. In Section 3 we show that any transverse family having the Desargues property admits a special family of symmetries, generalizing a certain $1$-parameter group of projective transformation (these symmetries are in fact {\em parabolic maps}, in the Dynamical terminology).    Finally, in Section 4 we show that symmetries of the Segre-Webster family obtained in Section 3 can not exist for a non-spherical hypersurface, and this proves \autoref{theor1}. Importantly, the symmetries we deal with here are {\em not} automorphisms of the hypersurface $M$ but are merely that of the Segre-Webster family, that is why one needs to develop here a rigidity theory {\em for transverse families of planar curves} (rather than real hypersurfaces). 

\bigskip

\begin{center} \bf Acknowledgments \end{center}

\bigskip

The author is grateful to Dmitri Zaitsev, Alexander Tumanov, and  Sergey Yakovenko for useful discussions, and also to Alexey Glyutsuk and Etien Gyss for their valuable comments on a talk of the author given in ENS Lyon in June 2015. 

 During the preparation of the paper the author was supported the Austrian Science Fund (FWF).

\bigskip

\section{Background material: Segre-Webster varieties}

Let $M$ be a smooth connected real-analytic hypersurface in
$\cx^{n+1}$, $Z=(z,w)\in \cx^n \times \cx$, $0\in M$, and let $U$ 
be a neighbourhood of the origin such that $M\cap U$ admits a real-analytic
defining function $\phi(Z,\overline Z)$ for $\phi$ holomorphic in $U\times \bar U$. For every point $\zeta\in
U$ we can associate to $M$ its so-called Segre-Webster variety in $U$
defined as
$$
Q_\zeta= \{Z\in U : \phi(Z,\overline \zeta)=0\}.
$$
Segre-Webster varieties depend holomorphically on the variable $\overline
\zeta$. One can find  a suitable pair of neighbourhoods $U_2={\
U_2^z}\times U_2^w\subset \cx^{n}\times \CC{}$ and $U_1 \Subset
U_2$ such that
$$
  Q_\zeta=\left \{(z,w)\in U^z_2 \times U^w_2: w = h(z,\overline \zeta)\right\}, \ \ \zeta\in U_1,
$$
is a closed complex analytic graph. Here $h$ is a holomorphic
function. Following \cite{DiPi} we call $U_1, U_2$ a {\it standard
pair of neighbourhoods} of the origin. The antiholomorphic
$(n+1)$-parameter family of complex hypersurfaces
$$\{Q_\zeta\}_{\zeta\in U_1}$$ is called \it the Segre-Webster family of $M$
at the origin. \rm From the definition and the reality condition
on the defining function the following basic properties of Segre-Webster
varieties follow (we assume $Z,\zeta\in U_1$ below):
\begin{equation}\label{e.svp}
  Z\in Q_\zeta \ \Leftrightarrow \ \zeta\in Q_Z,
\end{equation}
\begin{equation*}
  Z\in Q_Z \ \Leftrightarrow \ Z\in M,
\end{equation*}
\begin{equation*}
  \zeta\in M \Leftrightarrow \{Z\in U_1: Q_\zeta=Q_Z\}\subset M.
\end{equation*}
The fundamental role of Segre-Webster varieties for holomorphic mappings
is illuminated by their invariance property: if $f: U \to U'$ is a
holomorphic map sending a smooth real-analytic hypersurface
$M\subset U$ into another such hypersurface $M'\subset U'$, and $U$
is as above, then
$$
f(Z)=Z' \ \ \Longrightarrow \ \ f(Q_Z)\subset Q'_{Z'}.
$$ For
the proofs of these and other properties of Segre-Webster varieties see,
e.g., \cite{webster}, \cite{DiPi},  or
\cite{ber}.

In the particularly important case when $M$ is a \it real
hyperquadric, \rm i.e., when $$M=\left\{
[\zeta_0,\dots,\zeta_N]\in \cx\mathbb P^{N} : H(\zeta,\bar \zeta)
=0  \right\},$$ where $H(\zeta,\bar \zeta)$ is a nondegenerate
Hermitian form in $\CC{N+1}$ with $k+1$ positive and $l+1$
negative eigenvalues, $k+l=N-1,\,0\leq l\leq k\leq N-1$, the Segre-Webster
variety of a point $\zeta \in\CP{N}$ is the projective hyperplane
$$ Q_\zeta = \{\xi\in \cx\mathbb P^N: H(\xi,\bar\zeta)=0\}.$$ (In particular, the $(2N-1)$-dimensional sphere $S^{2N-1}\subset\CC{N}$ falls into the category of real hyperquadrics). The
Segre-Webster family $\{Q_\zeta,\,\zeta\in\CP{N}\}$  coincides in this
case with the space $(\CP{N})^*$ of all projective hyperplanes in
$\CP{N}$.

The space of Segre-Webster varieties $\{Q_Z : Z\in U_1\}$ can be
identified with a subset of $\cx^K$ for some $K>0$ in such a way
that the so-called \it Segre-Webster map \rm $\lambda : Z \to Q_Z$ is
holomorphic. For a Levi nondegenerate at a point
$p$ hypersurface $M$, its Segre-Webster map is one-to-one in a
neighbourhood of $p$. When $M$ contains a complex hypersurface
$X$, for any point $p\in X$ we have $Q_p = X$ and $Q_p\cap
X\neq\emptyset\Leftrightarrow p\in X$, so that the Segre-Webster map
$\lambda$ sends the entire $X$ to a unique point in $\CC{K}$ and,
accordingly, $\lambda$ is not even finite-to-one near each $p\in
X$ (i.e., $M$ is \it not essentially finite \rm at points $p\in
X$). For a hyperquadric $\Q\subset\CP{N}$ the Segre-Webster map $\lambda'$
is a global natural one-to-one correspondence between $\CP{N}$ and
the space $(\CP{N})^*$.

\section{Symmetries of families with the Desargues property}

As was explained in the Introduction, the necessity of the Desargues condition in \autoref{theor1} is obvious, that is why the rest of the paper is dedicated to the proof of sufficiency of the Desargues condition in \autoref{theor1}.

Throughout this section we consider a transverse family $\sss$ in a polydisc $U$ centered at $O$ such that all admissible configurations of seven points in $U$ for $\sss$ have the Desargues property. In the latter setting, we aim to prove that the family $\sss$ has a ``large'' local automorphism group. We specify that a local biholomorphism $F:\,(\CC{2},O)\mapsto  (\CC{2},O)$ is a {\em local automorphism of $\sss$}, if it transforms any curve in $\sss$ into another curve in $\sss$. It is convenient for our purposes   to  assume that both polydiscs $U,V$ are centered at the origin (so that $O$ is the origin), and the curve $L_0$ is given by 
$$L_0=\{w=0\},$$ where $(z,w)\in\CC{}\times\CC{}$ are the coordinates in $\CC{2}$. (Here we keep the notations from the Introduction). 
\begin{definition}
The curves of the family $\sss$ are called in what follows {\em $\sss$-lines}. 
\end{definition}

 We then construct a family of automorphisms of $\sss$, {\em preserving the distinguished $\sss$-line $L_0$}. The idea behind  that is to think of $L_0$ as the ``line at infinity'', so that automorphisms preserving $L_0$ become analogues of affine transformations. It is natural then to call two $\sss$-lines $L_1,L_2$ {\em parallel}, if $L_1\cap L_2$ is non-empty and lies in $L_0$. We write $L_1||L_2$ in the latter case. In this language, the Desargues property reads as follows:

\bigskip

{\em let each of the triples of points $(O,A,A'),(O,B,B'),(O,C,C')$ belong to some (its own) $\sss$-line, and assume that $AB||A'B',BC||B'C'$; then $AC||A'C'$.}

 \bigskip 
 
Using geometric intuition arising from the above reformulation, we define the desired automorphisms $\varphi$ of $\sss$. These are analogues of the linear-fractional transformations
\begin{equation}\label{scalings}z\mapsto \frac{z}{1-rw},\quad w\mapsto \frac{w}{1-rw},\quad r\neq 0,\end{equation}
and have the following properties:

\bigskip

\noindent (a)  $\varphi(B)=B$ for all $B\in L_0$; in particular, $\varphi(O)=O$. 

\bigskip

\noindent (b) $\varphi(OB)\subset OB$ for all $B\neq O$, if $OB$ exists.

\bigskip

\noindent (c) the $\sss$-line $\varphi(B)\varphi(C)$ is parallel to the $\sss$-line $BC$, if both $\sss$-lines exist. 

\bigskip

We start with choosing a domain $\mathcal D$ of the form
$$\mathcal D=\bigl\{|z|<r_0,\quad |w|<c_0|z|\bigr\}$$
for sufficiently small $r_0,c_0>0$. We then have the property that $\mathcal D\subset U$ and that for all $B\in\mathcal D$ the $\sss$-line $OB$ exists. Next, we fix some $A\in \mathcal D,\,A\notin L_0$ and  $A'\in OA\cap \mathcal D$ sufficiently enough to $O$.  For the set of points $B\in\mathcal D$, not lying in $OA$ and sufficiently close to the origin, we then define a map $\varphi=\varphi (A,A')$.  For that, we consider the $\sss$-line $AB$ (the latter exist since the   $\sss$-line $AO$ exists and $B$ is close to $O$), and take its intersection $P$ with the ``line at infinity'' $L_0$ (the intersection, again, exists because the intersection of $AO$ with $L_0$ exists, and $B$ is close to $O$). We then put:
$$\varphi(B):=B'.$$
The construction is illustrated by the picture below. We emphasize again that $\varphi$ is defined so far on $\mathcal D\setminus OA$ only. Clearly, $\varphi$ is holomorphic wherever it is defined (as follows from the fact that the family $\sss$ is holomorphic and from the implicit function theorem).

\includegraphics[scale=0.4]{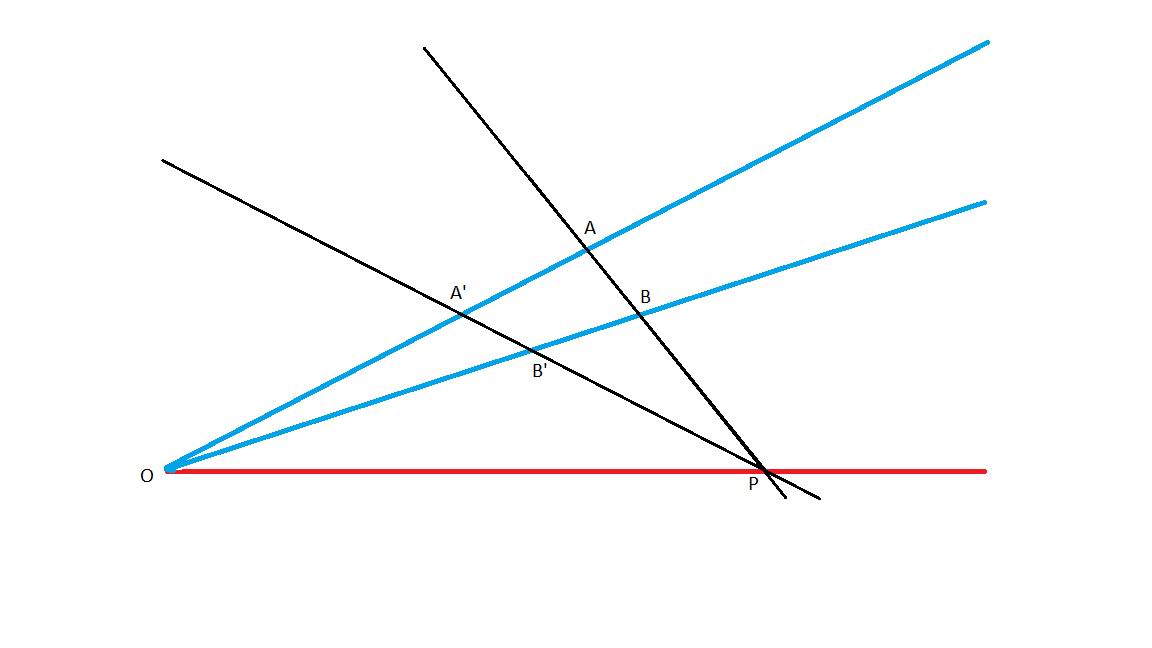}

 To extend $\varphi$ to $OA$ analytically, we take some $B\in\mathcal D\setminus OA$ and consider the map $\psi=\psi(B,B')$, defined in the domain $\mathcal D\setminus OB$   using the points $B,B'$ identically to $\varphi$.  Then, considering any point $C\in\mathcal D$ where both $\varphi,\psi$ are defined and putting $C':=\varphi(C),\,C'':=\psi(C)$, we get from the definition of $\varphi$ and the Desargues property:
$$ AB||A'B',\quad AC||A'C'\quad \Longrightarrow \quad BC||B'C'.$$
Similarly, from the definition of $\psi$ we get:
$$ BA||B'A',\quad BC||B'C''\quad \Longrightarrow \quad AC||A'C''.$$
Thus $$AC||A'C'||A'C'', \quad BC||B'C'||B'C'',$$
so that $$C''\in A'C', \quad C''\in B'C'.$$ We get $C''=C'$, and this proves that $\varphi$ and $\psi$ coincide wherever both are defined. Thus, $\psi$ is the desired analytic extension of $\varphi$ to $\mathcal D$, as required. It is immediate then (from the construction of $\varphi$ and the fact that $\varphi(A,A')=\psi(B,B')$ for any choice of $B,B'$)  that $\varphi$ satisfies (a),(b) and (c).  

The next step is to extend $\varphi$ holomorphically to a neighborhood of $O$. For that, we fix two distinct points $P,Q\in L_0$ both different from $O$, and fix two distinct $\sss$-lines $L_1=OX$ and $L_2=OY$ for some fixed $X,Y\in\mathcal D$ (the points $X,Y$ are actually of no interest to us). Then, for any point $K\in U\setminus \{L_0\cup L_1\cup L_2\}$ sufficiently close to $O$, consider the line $KP$ (the latter exists if $K$ is close enough to $O$). We then choose a point $B\in KP\cap L_1$, and set $B':=\varphi(B)$ (the intersection exists if $K$ is close enough to $O$). We have $B'\in OB$. Observe that, if $K\rightarrow O$, then  $B$ and hence $B'$ are also arbitrarily close to $O$, so that $PB'$ is non-empty (since $PO=L_0$ is non-empty and $P$ is fixed). We then repeat this construction by choosing  $C\in KQ\cap L_2$ and setting $C':=\varphi(C)$. We similarly have $QC'\neq\emptyset$. Note finally that the intersection of $PB'$ and $QC'$ is non-empty for $K$ sufficiently close to $O$ (since then $B'$ and $C'$ are arbitrarily close to $O$, while $P$ and $Q$ are fixed). Now we set $$K':=PB'\cap QC'.$$
It is then immediate to see from the Desargues property and the definition of $\varphi|_{\mathcal D}$ that for all $K\in\mathcal D\setminus \{L_0\cup L_1\cup L_2\}$ which are close enough to $O$, we have $\varphi(K)=K'$. Indeed, for such $K$ we have:
$$BK||B'K',\quad CK||C'K'.$$
Also, setting $K'':=\varphi(K)$, we have by the construction of $\varphi|_{\mathcal D}$:
$$ BK||B'K'',\quad CK||C'K''.$$
Hence $$K''\in B'K'=B'P,\quad K''\in C'K'=C'Q,$$ and we conclude that $K''=K'$, as required. The picture below illustrates the argument.

\includegraphics[scale=0.39]{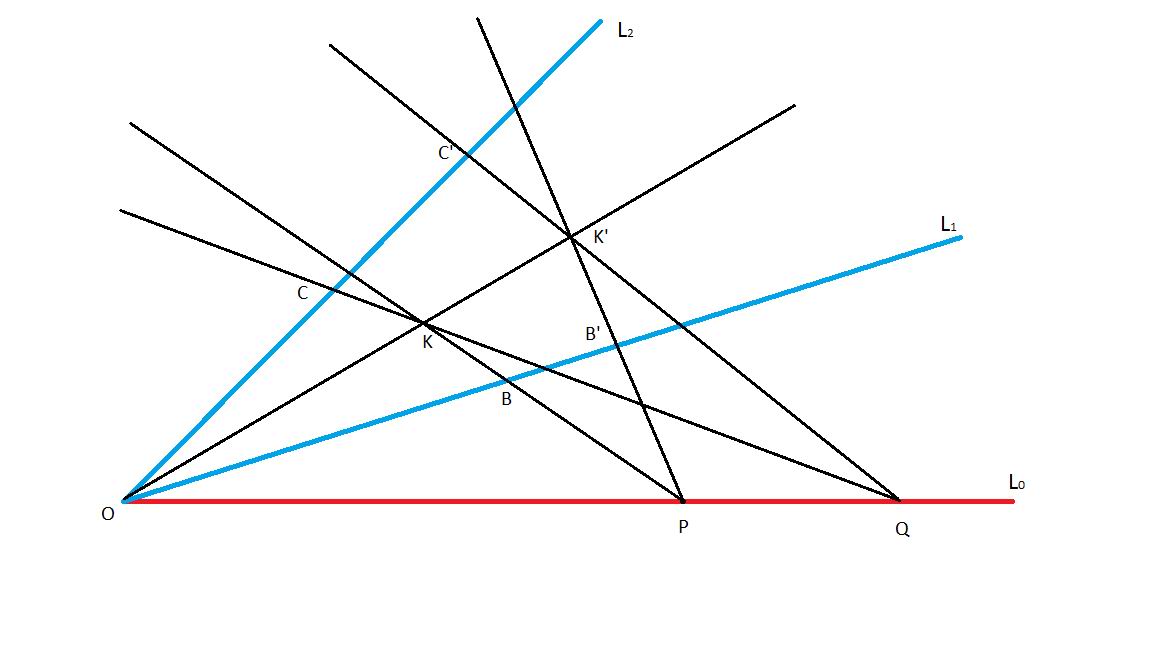}

Thus, setting 
$$\varphi(K):=K',$$
we get an analytic continuation of $\varphi$ from $\mathcal D$ to a sufficiently small neighborhood $V=V(\varphi)$ of $O$ in $\CC{2}$ (the holomorphic extension to $O$ itself is obtained by removing the compact singularity $\{O\}$). The fact that $\varphi$ sends $\sss$-lines to $\sss$-lines and satisfies the properties (a)-(c) can be verified  from the Desargues property of $\sss$ and the construction of $\varphi$ by arguments very similar to the ones above, so that we leave the details to the reader. By construction, $\varphi$ is one-to-one, thus it is an automorphism of the family $\sss$.

Let us now analyse the $1$-parameter family of automorphisms $\{\varphi (A,A')\}_{A'\in OA}$ constructed above. Since each $\varphi (A,A')$ transforms an $\sss$-line $OB,\,B\in\mathcal D$ into itself, we get that the differential of $\varphi (A,A')$ at $O$ has the form $$\lambda\cdot\mbox{Id}$$ for some $\lambda\neq 0$. On the other hand, $\varphi|_{L_0}$ is the identity in view of property (a), hence the differential of $\varphi$ at $O$   is the identity. Since all $\varphi (A,A')$ are distinct, and neither of them is the identity (for a fixed $A$ and various $A'\neq A$), we get a $1$-parameter family of local automorphisms of $\sss$, fixing $O$, different from the identity, having the identity linear part at $O$, and fixing each point in $L_0$. From here we obtain the following crucial statement, summarizing this section.

\begin{proposition}\label{key}
For any transverse holomorphic family $\sss$ in a polydisc $U$ centered at a point $O$ with a distinguished $\sss$-line $L_0\ni O$, there exists a $1$-parameter family of local automorphisms $\{\varphi_{r}\}_{r\in\Omega}$ of $\sss$, all preserving $O$, different from the identity, having the identity linear part at $O$, and being the identity when restricted onto $L_0$. Here $\Omega$ is a domain in $\CC{}$, and each $\varphi_{r}$ is defined in its own open neighborhood $V_r$ of $O$. 
\end{proposition}


\section{Proof of the main result}

In this section we apply \autoref{key} for proving \autoref{theor1}. An immediate intent here would be to apply one of the existing rigidity theorems in CR-geometry  saying that a Levi-nondegenerate hypersurface in $\CC{2}$, which possesses an automorphism $\varphi\neq\mbox{Id}$ with the identity linear part at some fixed point $O\in M$, is automatically spherical (e.g.,  \cite{belold,KrLo}). However, such an argument is not possible here since an automorphism of the Segre-Webster family $\sss$ of $M$ does {\em not} need  be an automorphism of $M$ (see, e.g., \cite{faran},\cite{divergence},\cite{nonanalytic} for considerations related to the latter fact). That is why a similar statement for automorphism of the Segre-Webster family must be established separately.   

\medskip

We start by bringing the initial real-analytic Levi-nondegenerate hypersurface $M\subset\CC{2}$ to the Chern-Moser normal form \cite{chern}. Recall that the latter has the form
\begin{equation}\label{nf}
v=z\bar z+\sum_{k,l\geq 2}\Phi_{kl}(u)z^k\bar z^l, \quad (z,w)=(z,u+iv)\in\CC{2},
\end{equation}
where the ``resonant''\, terms $\Phi_{kl},\,k,l\geq 2$ satisfy, in addition,
$$\Phi_{22}=\Phi_{23}=\Phi_{32}=\Phi_{33}=0.$$ 
We then transfer to the
so-called \it complex defining
 equation \rm (see, e.g., \cite{ber})\,
$$w=\Theta(z,\bar z,\bar w)$$  of $M$ near the origin, which is
obtained by substituting $u=\frac{1}{2}(w+\bar
w),\,v=\frac{1}{2i}(w-\bar w)$ into the real defining equation and
applying the holomorphic implicit function theorem. In terms of the complex defining equation, the normal form \eqref{nf} reads as:
\begin{equation}\label{nfc}
w=\bar w+2iz\bar z+\sum_{k,l\geq 2}\Theta_{kl}(\bar w)z^k\bar z^l, \quad \Theta_{22}=\Theta_{23}=\Theta_{32}=\Theta_{33}=0.
\end{equation}
Since the  Segre-Webster
variety $Q_p$ of a point $p=(\xi,\eta)$ is given in terms of $\Theta$
by
$$w=\rho(z,\bar \xi,\bar \eta),$$
we 
can (after an anti-holomorphic reparameterization ) use \eqref{nfc} to write down the Segre-Webster family $\sss$ of $M$ as 
\begin{equation}\label{segre}
w=b+az+\sum_{k,l\geq 2}\Psi_{kl}(b)z^k a^l, \quad \Psi_{22}=\Psi_{23}=\Psi_{32}=\Psi_{33}=0,
\end{equation}
where $a,b$ are holomorphic parameters near the origin
($\Psi_{kl}$ are obtained from $\Theta_{kl}$ by a simple linear transformation formula). 

Next, associated with \eqref{segre} is the submanifold $$\mathcal M\subset\CC{4}=\CC{2}\times\CC{2}$$ given by the same equation \eqref{segre} (where the coordinates in $\CC{4}$ are $(z,w,a,b)$). Note that $\mathcal M$ is linearly equivalent to the complexification $M^{\CC{}}$ of $M$. We further observe that any automorphism $$F(z,w)=\bigl(f(z,w),g(z,w)\bigr)$$  of the Segre-Webster family $\sss$ preserving the origin and the $\sss$-line $L_0=\{w=0\}$ generates a (unique!) similar transformation $$G(a,b)=\bigl(\lambda(a,b),\mu(a,b)\bigr)$$ in the space of parameters such that the resulting product transformation $$\bigl(F(z,w),G(a,b)\bigr):\,(\CC{4},0)\mapsto(\CC{4},0)$$ preserves $\mathcal M$. (To obtain $G(a,b)$ one needs to simply take the unique $(a',b')$ such that the germs at the origin of the $\sss$-lines $F(L_{a,b})$ and $L_{a',b'}$ coincide).

We apply the latter observation to some automorphism $F=\bigl(f(z,w),g(z,w)\bigr)$ arising from \autoref{key}, where we choose the distinguished line to be 
$$L_0:=\{w=0\}.$$
We have 
\begin{equation}\label{specialmap}
dF|_0=\mbox{Id},\quad F(z,0)=(z,0),\quad F\neq\mbox{Id}.
\end{equation}
From the above argument, we get a biholomorphism $G=\bigl(\lambda(a,b),\mu(a,b)\bigr)$ of $(\CC{2},0)$ such that the product biholomorphism $\bigl(F(z,w),G(a,b)\bigr):\,(\CC{4},0)\mapsto(\CC{4},0)$ preserves the manifold $\mathcal M$, given by \eqref{segre}.  

We now assume that $M$ is non-spherical, so that at least one $\Psi_{kl}$ in \eqref{segre} is not identically zero. Then we  perform a power series calculation giving a contradiction  with the existence of the local automorphism $\bigl(F(z,w),G(a,b)\bigr)$ of $\mathcal M$. For that, let us introduce for the coordinates in the product space $\CC{4}=\CC{2}\times\CC{2}$ the weights:
\begin{equation}\label{weights}
[z]=[a]=1,\quad [w]=[b]=2.
\end{equation}
Then the components of the product map $$\bigl(F(z,w),G(a,b)\bigr)=\bigl(f(z,w),g(z,w),\lambda(a,b),\mu(a,b)\bigr)$$ admit the expansion:
\begin{equation}\label{expansion}
\begin{aligned} 
f(z,w)=z+\sum_{j\geq 2}f_j(z,w),&\quad g(z,w)=w+\sum_{j\geq 3}g_j(z,w),\\ 
\lambda(a,b)=\sum_{j\geq 1}\lambda_j(a,b),&\quad \mu(a,b)=\sum_{j\geq 1}\mu_j(a,b),
\end{aligned}
\end{equation}
where $f_j,g_j,\lambda_j,\mu_j$ are homogeneous polynomials of a fixed weight $j$ with respect to the gradation \eqref{weights} (the special form of $f,g$ follows from  \eqref{specialmap}). Similarly, we expand $\mathcal M$ as:
$$w=b+az+\sum_{j\geq 6}\Psi_j(z,a,b)$$
($\Psi_j$ here is a homogeneous polynomial in $z,a,b$ of the uniform weight $j$).
Next, let us denote by $m$ the smallest integer $j\geq 6$ such that $\Psi_j\not\equiv 0$ ($m$ is finite by the non-sphericity assumption). Thus we can write $\mathcal M$ as:
\begin{equation}\label{manifold}
w=b+az+\sum_{j\geq m}\Psi_j(z,a,b),
\end{equation}
where $$\Psi_m(z,a,b)\neq 0.$$

\medskip

The fact that $\bigl(F(z,w),G(a,b)\bigr)$ preserves $\mathcal M$ reads as:
\begin{equation}\label{identity}
g(z,w)=\Psi\bigl(f(z,w),\lambda(a,b),\mu(a,b)\bigr)|_{w=\Psi(z,a,b)}
\end{equation}
(here $\Psi(z,a,b)$ is the right-hand side of \eqref{manifold}). 
We now substitute \eqref{expansion},\eqref{manifold} into \eqref{identity}, and for each fixed weight $j\geq 1$ collect all terms of weight $j$. Then we use the $j$-th identity to compute the collection
$$\bigl(f_{j-1}(z,w),g_j(z,w),\lambda_{j-1}(a,b),\mu_j(a,b)\bigr).$$  For $j=1$ we obtain:  
\begin{equation}\label{j1}\mu_1=0\end{equation}
For $j=2$ we obtain:  
\begin{equation}\label{j2}\lambda_1=a,\quad \mu_2=b.\end{equation}
For $j=3$ we obtain:  
\begin{equation}\label{j3}f_2=0,\quad g_3=0,\quad \lambda_2=0,\quad \mu_3=0\end{equation}
(to get $f_2=0$ we  used \eqref{specialmap}).
For $j=4$ we obtain:  
\begin{equation}\label{j4}
f_3=rzw,\quad g_4=rw^2,\quad \lambda_2=rab,\quad \mu_3=rb^2, \quad r\in\CC{}
\end{equation}
(we again used \eqref{specialmap}). We now summarize the above calculations, taking \eqref{specialmap} into account:
\begin{equation}\label{summary}
\begin{aligned}
f_1=z,\quad \lambda_1=a,\quad g_1=\mu_1=f_2=\lambda_2=g_3=\mu_3=0,\\
f_3=rzw,\quad g_4=rw^2,\quad \lambda_2=rab,\quad \mu_3=rb^2.
\end{aligned}
\end{equation}

Thus, we conclude that, importantly, 

\medskip

{\em the map $F(z,w)$ is approximated by an appropriate map \eqref{scalings}} 

\medskip

(where $f$ is approximated up to weight $3$, and $g$ up to weight $4$).

\medskip

We claim now that we have $r\neq 0$ in \eqref{summary}. To prove that, we provide some homological argument (in the spirit of Poincar\'e-Moser). Assume that, otherwise, $r=0$. Let us consider a weight $l\geq 1$ identity arising from \eqref{identity}, and extract from it all terms  involving the collection $(f_{l-1},g_l,\lambda_{l-1},\mu_l)$. It is not difficult to check that we get:
\begin{equation}\label{weightl}
\mathcal L(f_{l-1},g_{l},\lambda_{l-1},\mu_{l})=T_{l}\bigl(\{f_{j-1},g_{j},\lambda_{j-1},\mu_{j}\}_{j<l}\bigr).
\end{equation}
We explain the notations in \eqref{weightl}:
here $\mathcal L$ is the linear operator
\begin{equation}\label{L}
\mathcal L(f,g,\lambda,\mu):=g(z,b+az)-\mu(a,b)-af(z,b+az)-z\lambda(a,b),
\end{equation}
and $$T_l\bigl(\{f_{j-1},g_{j},\lambda_{j-1},\mu_{j}\}_{j<l}\bigr)$$ is some precise polynomial of 
$\{f_{j-1},g_{j},\lambda_{j-1},\mu_{j}\}_{j<l}$  and their derivatives, exact form of which is of no interest to us (this form in fact depends on a concrete $\Psi$).

Let us now introduce the linear  space $\mathcal U$ of collections $(f,g,\lambda,\mu)$ satisfying \eqref{summary} {\em with $r=0$}. We then  consider the image $\mathcal V$ of the operator $\mathcal L$  restricted onto the space $\mathcal U$. Denote also by $\mathcal W$ the space of power series of the form $$\sum_{j\geq 3} \Psi_j(z,a,b)$$ (we still use here the gradation \eqref{weights}). Then, arguing very similarly to the proof of Lemma 2.1 in \cite{chern}, we can prove the following
\begin{lemma}\label{normal}
\mbox{}

(i) The linear operator $\mathcal L$ is injective on $\mathcal U$.

\smallskip

(ii) The space $\mathcal W$ can be decomposed as
\begin{equation}\label{direct}
\mathcal W=\mathcal V\oplus \mathcal N,
\end{equation}
where $\mathcal N$ for  power series $\Psi(z,a,b)$ expanded as $$\Psi(z,a,b)=\sum_{k,l\geq 0}\Psi_{kl}(b)z^ka^l$$ is given by the conditions:
$$\Psi_{k0}=\Psi_{0l}=\Psi_{k1}=\Psi_{1l}=\Psi_{22}=\Psi_{23}=\Psi_{32}=\Psi_{33}=0,\quad k,l\geq 0.$$
\end{lemma}

Now, combining \eqref{weightl}, claim (i) in \autoref{normal}, and the assumption $r=0$ for the automorphism $(f,g,\lambda,\mu)$ under consideration, we immediately conclude that {\em such an automorphism with $r=0$ is unique}. Indeed, in each step $l\geq 1$ we solve a non-homogeneous system of linear equations for the step-$l$-collection $(f_{l-1},g_l,\lambda_{l-1},\mu_l)$, while statement (i) in \autoref{normal} says that the corresponding homogeneous linear system has only the zero solution. Thus, our series of weighted identities  determines an automorphsim uniquely.  

We conclude from here our automorphism is the identity, which is a contradiction with \eqref{specialmap}. This proves finally that $$r\neq 0\quad\mbox{in}\quad\eqref{summary}.$$

To compute $(f_{j-1},g_j,\lambda_{j-1},\mu_j)$ with $5\leq j\leq m-1$ we argue as follows. It is a well-know fact (which, in fact, can also be verified from \autoref{normal}) that any (local) map preserving the family of straight lines in $\CC{2}$ is linear-fractional (i.e., it is a restriction of an automorphism of $\mathbb{CP}^2$; see, e.g., \cite{tresse}). Any such map, preserving the origin and satisfying \eqref{specialmap}, is clearly given by \eqref{scalings}.  At the same time, it is easy to check that the product map
\begin{equation}\label{g+}
z\mapsto  \frac{z}{1-rw}\quad w\mapsto \frac{w}{1-rw}, \quad a\mapsto \frac{a}{1-r b},\quad b\mapsto \frac{b}{1-rb}
\end{equation} 
preserves the ``model''\, manifold
\begin{equation}\label{model}
\mathcal M_0=\{w=b+az\},
\end{equation}
corresponding to the family of straight lines in $\CC{2}$. Thus, 

\bigskip

{\em all local automorphisms of $\mathcal M_0$ satisfying \eqref{specialmap} are given by \eqref{g+}}.

 \bigskip 
 
We continue the proof, still considering the initial automorphism of $\mathcal M$ satisfying \eqref{summary} with some $r\neq 0$, as well as an automorphism \eqref{scalings} of the model $\mathcal M_0$ {\em corresponding to the same value of $r$}. Note that, since $\mathcal M$ and $\mathcal M_0$ coincide up to weight $m-1$, the weight $j$ equations to compute $(f_{j-1},g_j,\lambda_{j-1},\mu_j)$ are the same for the model and the non-model cases respectively, and we conclude that each collection $(f_{j-1},g_j,\lambda_{j-1},\mu_j)$ with $j\leq m-1$ simply {\em coincides with that for the map \eqref{scalings}} with the same value of $r$. This gives the desired answer for $(f_{j-1},g_j,\lambda_{j-1},\mu_j)$ with $5\leq j\leq m-1$. 

Next, for $j=m$ and $j=m+1$ is is straightforward to check that the equations to obtain $(f_{j-1},g_j,\lambda_{j-1},\mu_j)$ are identical to that for the case of $\mathcal M_0$, so that the corresponding collections of weighted polynomials again {\em coincide} (the resonant terms $\Psi_m,\Psi_{m+1}$ cancel each other in the $m$-th and $(m+1)$-st identity respectively). 

Finally, we consider the weight $j=m+2$. Let us consider first the weight $m+2$ identity for the model $\mathcal M_0$. It has the form:
\begin{equation}\label{right}
\mathcal L(f_{m+1},g_{m+2},\lambda_{m+1},\mu_{m+2})=H_{m+2}\bigl(\{f_{j-1},g_{j},\lambda_{j-1},\mu_{j}\}_{j<m+2}\bigr),
\end{equation}
where $\mathcal L$ is the above linear operator \eqref{L},
and $$H_{m+2}\bigl(\{f_{j-1},g_{j},\lambda_{j-1},\mu_{j}\}_{j<m+2}\bigr)$$ is some precise polynomial of 
$\{f_{j-1},g_{j},\lambda_{j-1},\mu_{j}\}_{j<m+2}$ and their derivatives, exact form of which is of no interest to us.
Now, returning to the case of the non-model manifold $\mathcal M$ and considering the weight $m+2$ in \eqref{identity}, we obtain:
\begin{equation}\label{wrong}
\mathcal L(f_{m+1},g_{m+2},\lambda_{m+1},\mu_{m+2})=H_{m+2}\bigl(\{f_{j-1},g_{j},\lambda_{j-1},\mu_{j}\}_{j<m+2}\bigr)-rza\Psi_{m}
\end{equation}
(in this step the cancellation of resonant terms does not happen thanks to the non-linear terms in $f,g$ !).
Importantly, 

\medskip

{\em the expressions $H_{m+2}\bigl(\{f_{j-1},g_{j},\lambda_{j-1},\mu_{j}\}_{j<m+2}\bigr)$ corresponding to \eqref{right} and \eqref{wrong} coincide}

 \medskip
 
  \noindent (since, as was explained above, the collections $\{f_{j-1},g_{j},\lambda_{j-1},\mu_{j}\}_{j<m+2}$ coincide for the model and the non-model cases respectively).   
  
  To show now that the identities \eqref{right},\eqref{wrong} can not hold simultaneously, we make a simple observation that the difference of the two collections $(f,g,\lambda,\mu)$ (having the same $r\neq 0$) corresponding to $\mathcal M$ and $\mathcal M_0$ respectively belongs to the above space $\mathcal U$ (since this difference satisfies \eqref{summary} with $r=0$). Hence, subtracting \eqref{wrong} from \eqref{right} and using the notations of \autoref{normal}, we get an identity in which the left-hand side is in $\mathcal V$, while the right hand side is in $\mathcal N$. (The latter follows from the fact that, as the explicit description of $\mathcal N$ and the normal form \eqref{segre} show, we have the inclusions $$\Psi_{m}\in\mathcal N\quad \mbox{and}\quad za\Psi_m\in\mathcal N).$$ This, in  view of \eqref{direct} and the fact that $r\neq 0$, implies
$$\Psi_m(z,a,b)\equiv 0$$
and gives a contradiction.

\autoref{theor1} is proved now.

\end{document}